\documentclass[a4paper]{amsart}
\usepackage{amsmath}

\usepackage{amssymb}
\usepackage[ a4paper,left=15mm,top=15mm,bottom=15mm, right=15mm]{geometry}
\usepackage[all]{xy}
\usepackage{amsmath, amsfonts, amsthm , amssymb, amscd}
\usepackage{hyperref}
\usepackage{booktabs}
\usepackage{todonotes}
\usepackage{mathtools}
\usepackage{enumitem}
\usepackage{romannum}
\usepackage{anyfontsize}
\usepackage{tikz}
\usetikzlibrary{decorations.pathmorphing}
\newtheorem*{theorem*}{Theorem A}
\newtheorem*{theorem**}{Theorem B}

\newtheorem{theorem}{Theorem}[section]

\newtheorem{proposition}[theorem]{Proposition}
\newtheorem{lemma}[theorem]{Lemma}

\theoremstyle{definition}
\newtheorem{definition}[theorem]{Definition}

\newtheorem{problem}[theorem]{Problem}
\theoremstyle{remark}

\title{Orthogonality  related to Steiner and Pappus chains}

\author{Djordje Barali\'{c}, Vladimir Bo\v{z}ovi\'{c} \and Nikola Radoji\v{c}i\'{c}}
\address{\scriptsize{ Mathematical Institute SASA, Knez Mihailova 36, p.p. 367, 11001, Belgrade, Serbia }}
\email{djbaralic@mi.sanu.ac.rs}
\address{\scriptsize{Faculty of Science, University of Montenegro, P.O. Box 39040, Podgorica,  Montenegro}}
\email{vladobozovic@gmail.com}
\address{\scriptsize{Faculty of Science, University of Montenegro, P.O. Box 39040, Podgorica,  Montenegro}}
\email{nikola.radojicic@os-mnikcevic.edu.me }


\date{}

\begin{document}

\pagenumbering{arabic}

\begin{abstract}

The focus of this paper is on the study of specific circle formations known as orthogonal Pappus chains and the related incidence results that involve points of tangency between the circles in the construction. These chains give rise to new circle families, with their centres on a conic. The paper also explores similar questions for Steiner chains, although it is demonstrated that orthogonality cannot be defined for them. Nonetheless, explicit examples are provided, exhibiting comparable incidence results and raising whether two Steiner chains that share a common circle in the chains possess this property.

\end{abstract}

\maketitle

\section{Introduction}

\begin{figure}[!ht]
  \centering
    \includegraphics[width=0.4 \textwidth]{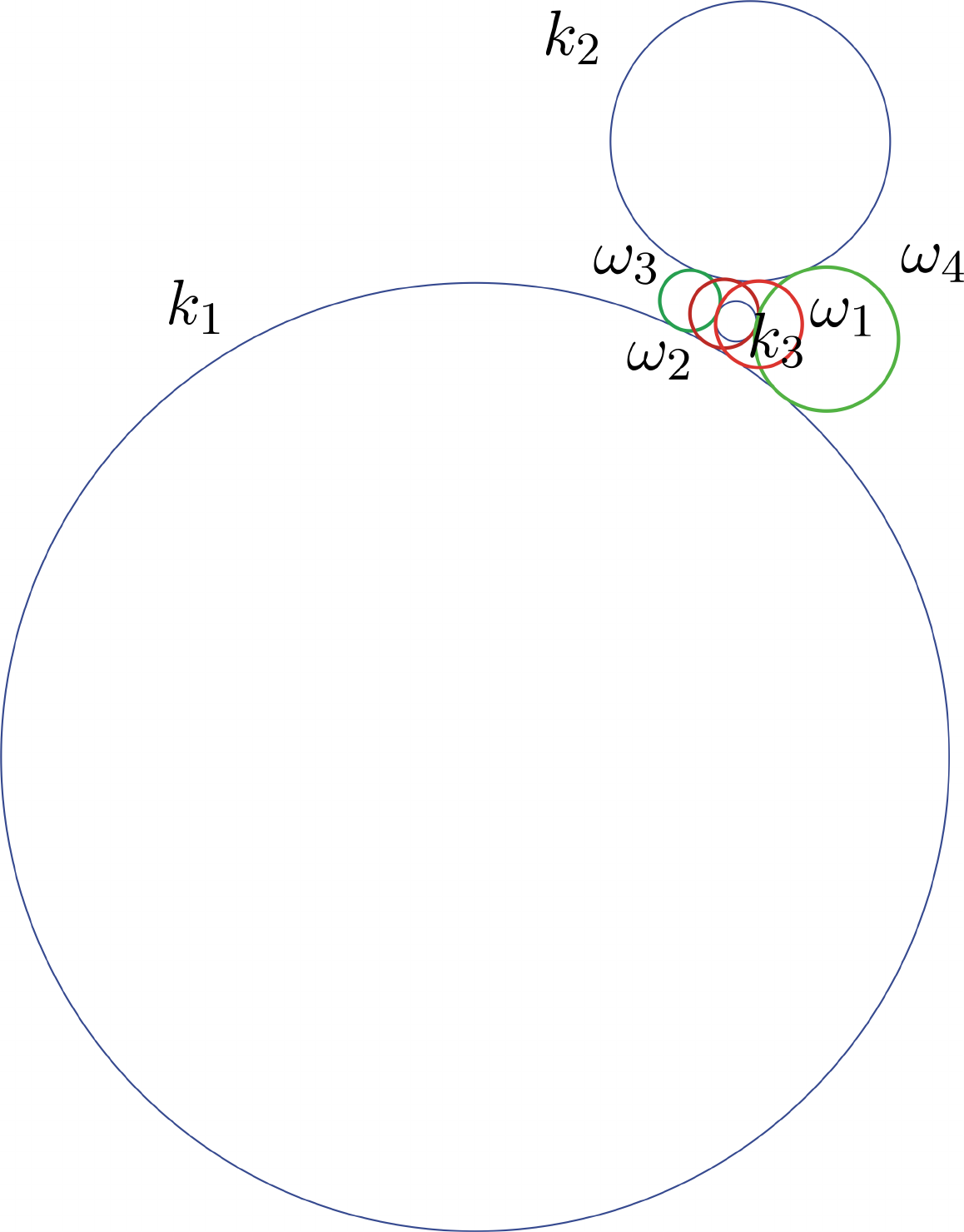}
    \caption{10th Apollonius' Problem}\label{Apol}
\end{figure}

Studying configurations of mutually tangent circles and conics has a long and exciting history in mathematics. The famous tenth problem of Apollonius of Perga (240-190) was the construction of four circles tangent to three given circles, Figure \ref{Apol}. Apollonius gave geometrical definitions of non-degenerate conics: ellipse, hyperbola and parabola, which we use today. Fascinating result in enumerative geometry on 3264 conics tangent to five given conics was posted by the famous geometer Jacob Steiner \cite{St1} and proved in 1859 by Jonqui\`{e}res, who did not publish it and independently in 1864 by Chasles in \cite{Chasles}. However, the first complete proof without gaps in classical intersection theory was given one century later by Fulton and MacPherson in \cite{Fulton}. The tenth Apollonius's problem and Steiner's conics problem ask for similar geometrical constraints, but the natural setting for the first is inversion and enumerative algebraic geometry for the second. In the 3rd century,  Pappus of Alexandria considered a ring of circles between two tangent circles,  today known as the Pappus chain,  such that each circle in the chain is tangent to the previous and the next circle in the chain. The prominent result on the Pappus chain, such that the centres two tangent circles and the first circle in the chain are collinear, states that the distance of the centre of the $n$-th circle in the chain from the line through the centre two circles is $n$ times of its radius see \cite{Pedoe}. Closely related to the Pappus chain is the Steiner chain of circles, where two circles tangent to all circles in the chain are not tangent but non-intersecting. Steiner himself discovered the renowned result on the Steiner chain in \cite{St},  proving that if there is a Steiner chain of $n$ circles such that the last circle in the chain is tangent to the first, then any circle tangent to the two given circles is a member of a Steiner chain of mutually tangent $n$ circles.

Steiner's theorem is an example of a poristic result. A geometrical configuration of figures is said to have a poristic property if it cannot exist unless a condition is met but otherwise may exist in an infinite number. Porisms have been considered one of the most profound results in classical Euclidean geometry. Celebrated Poncelet's porism is one of the most notable achievements in projective geometry. At the same time, its generalizations and variations have been extensively studied in the last decades due to their remarkable role in the theory of billiards, mechanics, and integrable dynamical systems, see \cite{DraRad} and \cite{DraRad1}. Several other
closure theorems are presented in \cite[Lecture~29]{Fuchs}.

In the recent paper \cite{Viera} \v{C}er\v{n}anov\'{a} studied the configuration called a fabric of kissing circles, obtained after applying circle inversion on a square grid filled with circles. She got riveting results about the curvatures of the circles in a fabric. In \cite{Lu} and \cite{Bar}, some further incidence properties of the Steiner chain configurations were studied. These results motivated us to investigate orthogonal Pappus chains and to introduce a more general configuration of orthogonal Pappus chains. Circle inversion is again our primary method for proving the claims.

\section{Pappus and Steiner chains}

Let $l$ and $m$ be two circles in a plane. Let $k_1$, $k_2, \dots, k_n, \dots $ be a sequence of circles tangent to $l$ and $m$ such that circles $k_i$ and $k_{i+1}$ are mutually tangent for all $i \geq 0$.

\begin{definition} \begin{itemize}
                     \item If circles $k$ and $l$ are mutually tangent, the sequence $k_1$, $k_2, \dots, k_n, \dots $ is called a Pappus chain.
                     \item If circles $k$ and $l$ do not intersect, the sequence $k_1$, $k_2, \dots, k_n, \dots $ is called a Steiner chain.
                   \end{itemize}
\end{definition}

In \cite{Lu} and \cite{Bar}, it was shown that certain interesting incidence results among points, lines and circles hold for Steiner chains. For all $i$, we label the following points:
\begin{itemize}
  \item By $L_i$ the tangent point of $k_i$ and $l$;
  \item By $M_i$ the tangent point of $k_i$ and $m$;
  \item By $N_i$ the tangent point of $k_i$ and $k_{i+1}$;
  \item By $t_i$ the common tangent of $k_i$ and $k_{i+1}$ at $N_i$;
  \item By $p_i$ the line $L_i M_i$.
\end{itemize}

In \cite[Lemma~3.1]{Bar}, it is proved that all lines $p_i$ and $t_i$ intersect in the same point that belongs to the line through the centres of circles $l$ and $m$ for the case of a Steiner chain. Here, we show that the same observation holds for a Pappus chain as well.

\begin{figure}[!ht]
  \centering
    \includegraphics[width= 0.7 \textwidth]{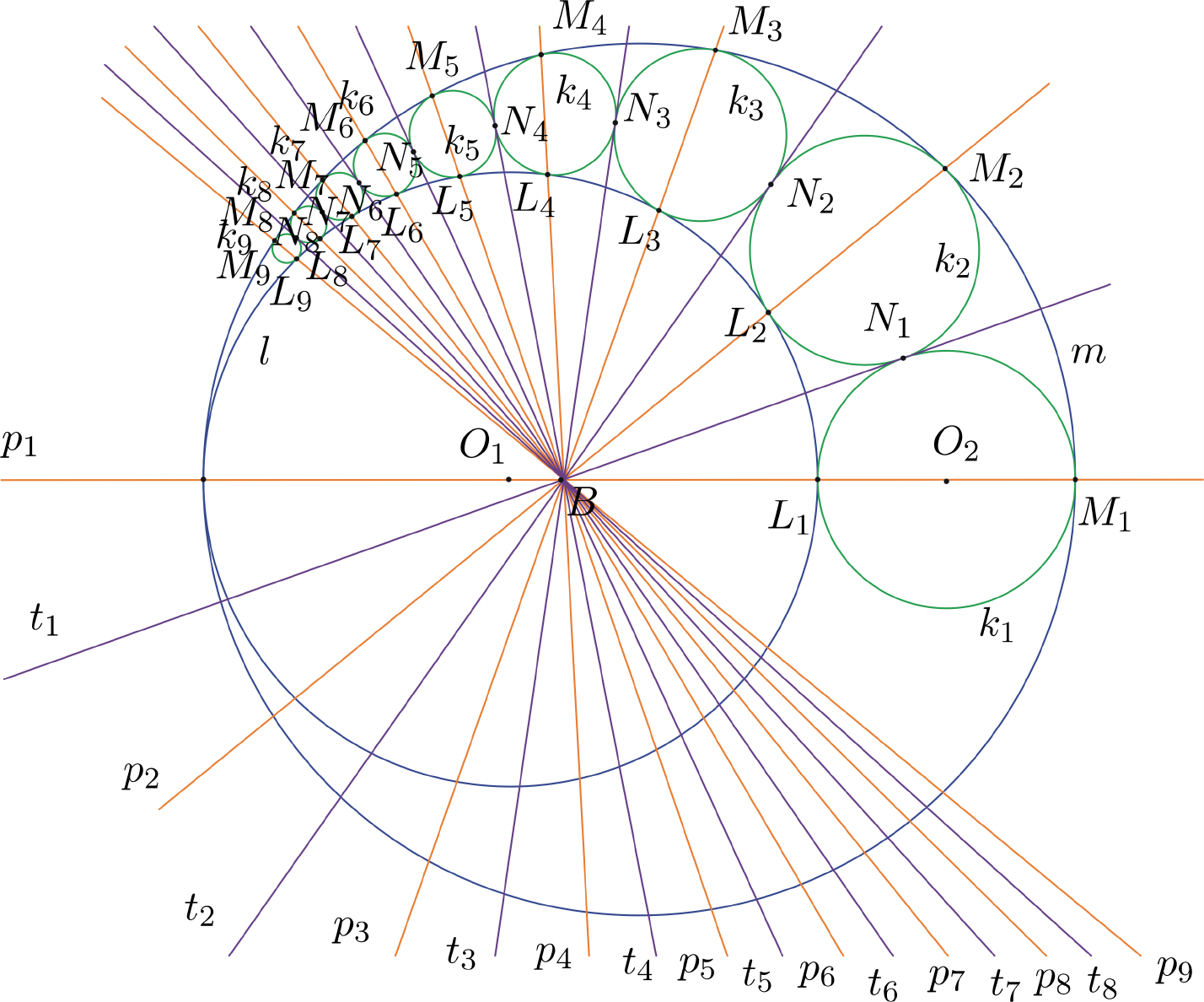}
    \caption{Theorem \ref{teo1}}\label{papus1}
\end{figure}

\begin{theorem}\label{teo1} For a Pappus chain, the lines $p_i$ and $t_i$ intersect in the same point $B$ that lies on the line through the  centers of $l$ and $m$, Figure \ref{papus1}.
\end{theorem}

\begin{figure}[!ht]
  \centering
    \includegraphics[width= 0.65 \textwidth]{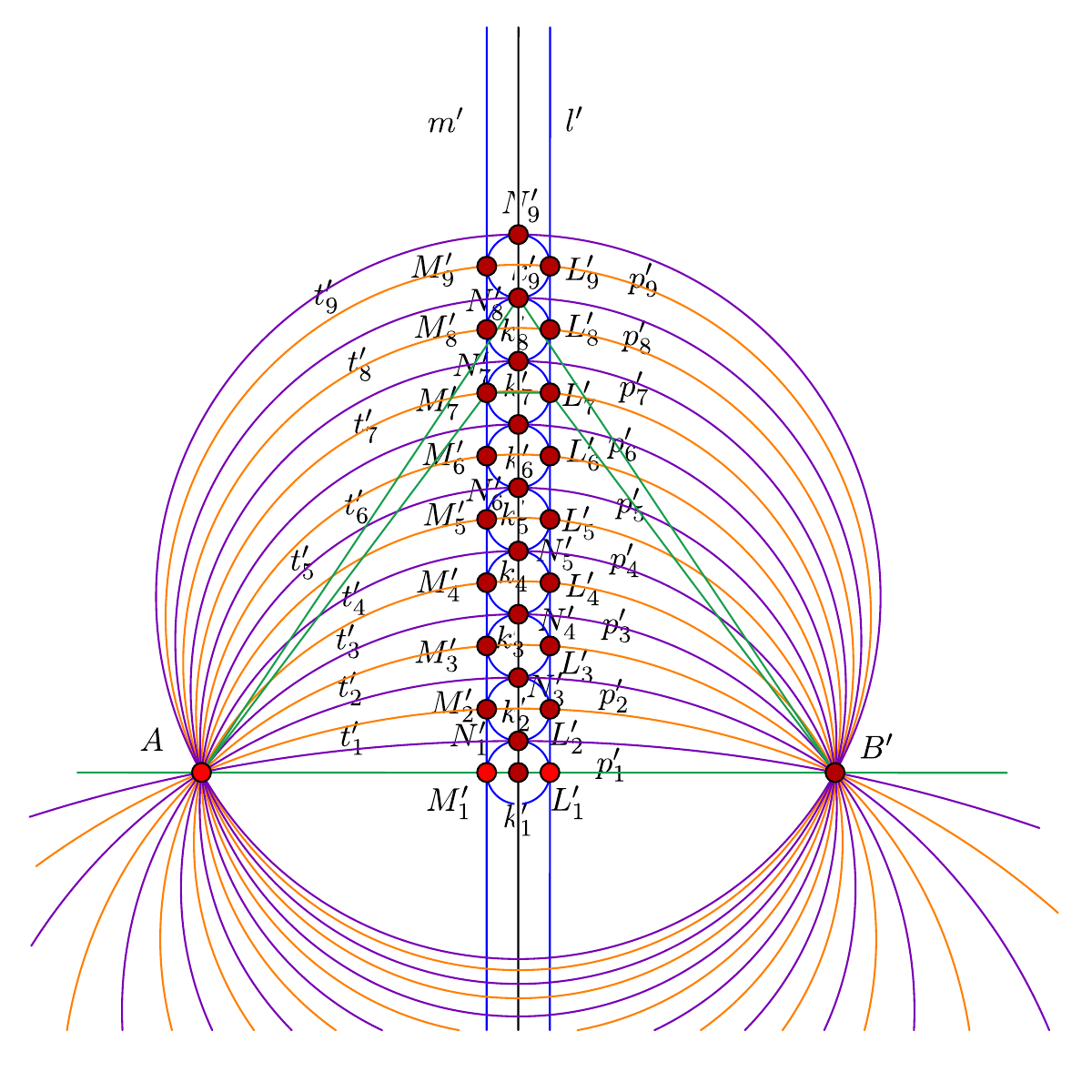}
    \caption{Proof of Theorem \ref{teo1}}\label{papus1i}
\end{figure}

\textsl{Proof:} Let the circles $l$ and $m$ touch at the point $A$. We apply an inversion with respect to a circle $\psi$ whose centre is $A$, the tangency point of $l$ and $m$,  Figure \ref{papus1i}. The circles $l$ and $m$ are mapped to parallel lines $l'$ and $m'$, orthogonal to the line $a$ passing through $A$ and the centres of $l$ and $m$. The chain of circles $k_i$ is mapped to the chain of congruent circles touching the lines $l'$ and $m'$. The lines $p_i$ are mapped to the circles $p'_i$ passing through $A$ for all $i$. Indeed, for any $i$, the line $L'_i M'_i$ is perpendicular to the lines $l'$ and $m'$, so $p'_i$ intersects the line in point $B_i$ such that the quadrilateral $A L'_i M'_i B_i$ is an isosceles trapezium. Therefore, the point $B_i$ is symmetric to $A$ concerning the perpendicular bisector of $L'_i M'_i$. But, the perpendicular bisector passes through the centres of all $k'_i$ and $N'_i$, so $B_i$ is fixed point $B$ and does not depend on $i$.

On the other hand, for all $i$, the line $t_i$ is mapped to the circle $t'_i$, tangent to $k'_i$ and $k'_{i+1}$ and passes through $A$. But, it means that the second intersection point $C_i$ of the line $a$ and $t'_i$ is such that $\triangle A N'_i C_i$ is isosceles. However, it is equivalent to $C_i$ being symmetrical to $A$ with respect to the line through the centres of the circles $k'_i$ in the chain, so $C_i$ is the same point $B$, and our proof is finished. \hfill $\square$

\begin{lemma}\label{lema}  For a Pappus chain, there exists a circle orthogonal to all circles $k_i$ in the chain, passing through all points $N_i$, and its centre is the common point of the lines $t_i$.
\end{lemma}

 \begin{figure}[!ht]
  \centering
    \includegraphics[width= 0.65 \textwidth]{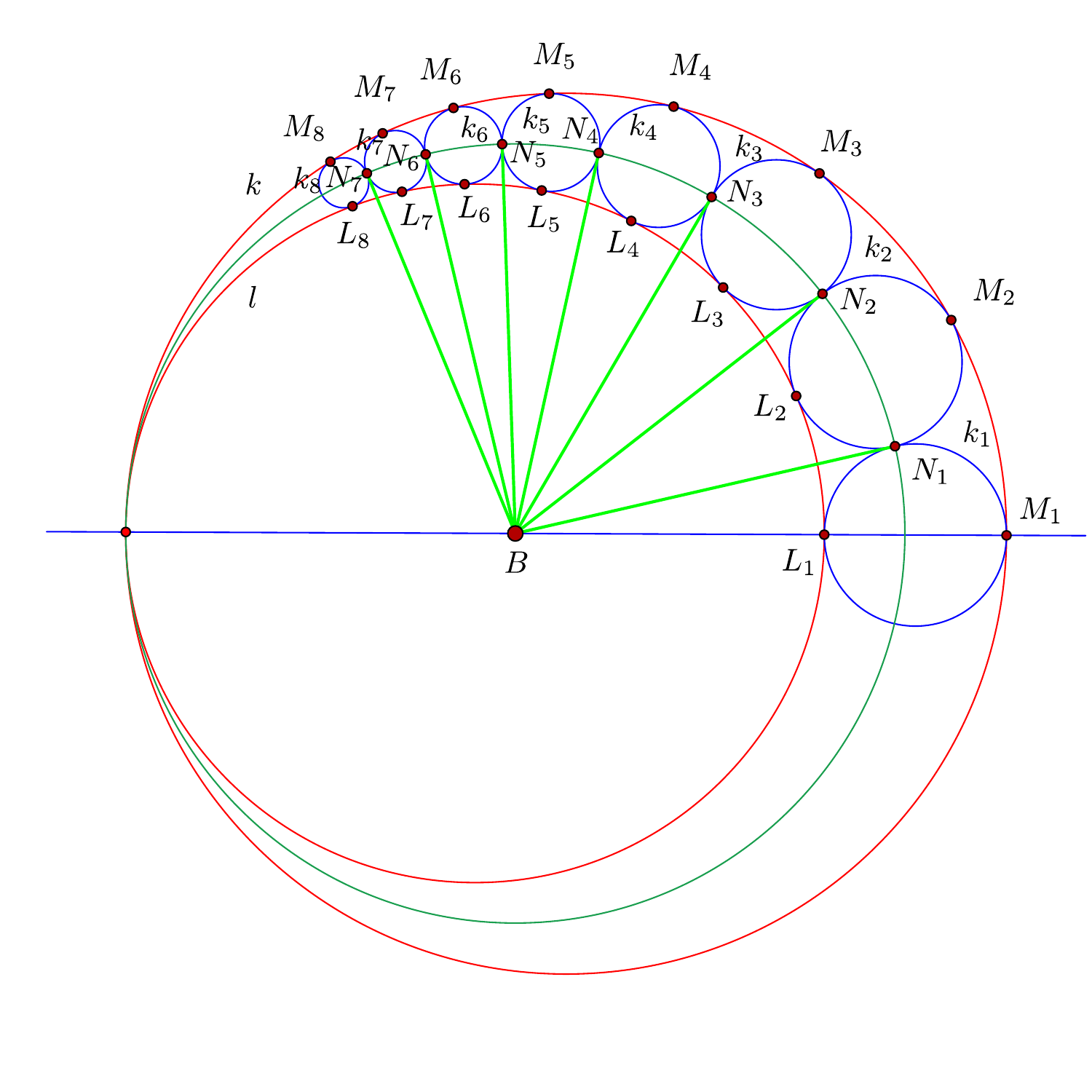}
    \caption{Lemma \ref{lema}}\label{figl}
\end{figure}

\textsl{Proof:} We apply the same inversion $\psi$ as in the previous proof. The line through the centres of all $k'_i$ is orthogonal to each circle in the chain and passes through all tangency points $N'_i$. Also, it is parallel to the lines $l'$ and $m'$ and orthogonal to the line $a$. Therefore, the preimage of this line is a circle tangent to the circles $l$ and $M$ at $A$, orthogonal to all circles $k_i$ in the chain and passing through all points $N_i$. From Theorem \ref{teo1}, it follows that all tangent segments $B N_i$ are equal, so $B$ is the centre of this circle, Figure \ref{figl}.   \hfill $\square$

\begin{proposition} \label{prop1} For all $j >i \geq 0$ the points $L_j, L_i, M_i, M_j$ lie on  the same circle $\varpi_{i, j}$ that meets circles $k_i$ and $k_j$ under angle $\varphi_{i,j}=\arctan (j-i)$.
\end{proposition}

\textsl{Proof:} Applying the inversion from the proof of Theorem \ref{teo1}, we easily deduce that the quadrilateral $L'_j L'_i  M'_i M'_j$ is a rectangle and its vertices lie on the same circle.
 \begin{figure}[!ht]
  \centering
    \includegraphics[width= 0.6 \textwidth]{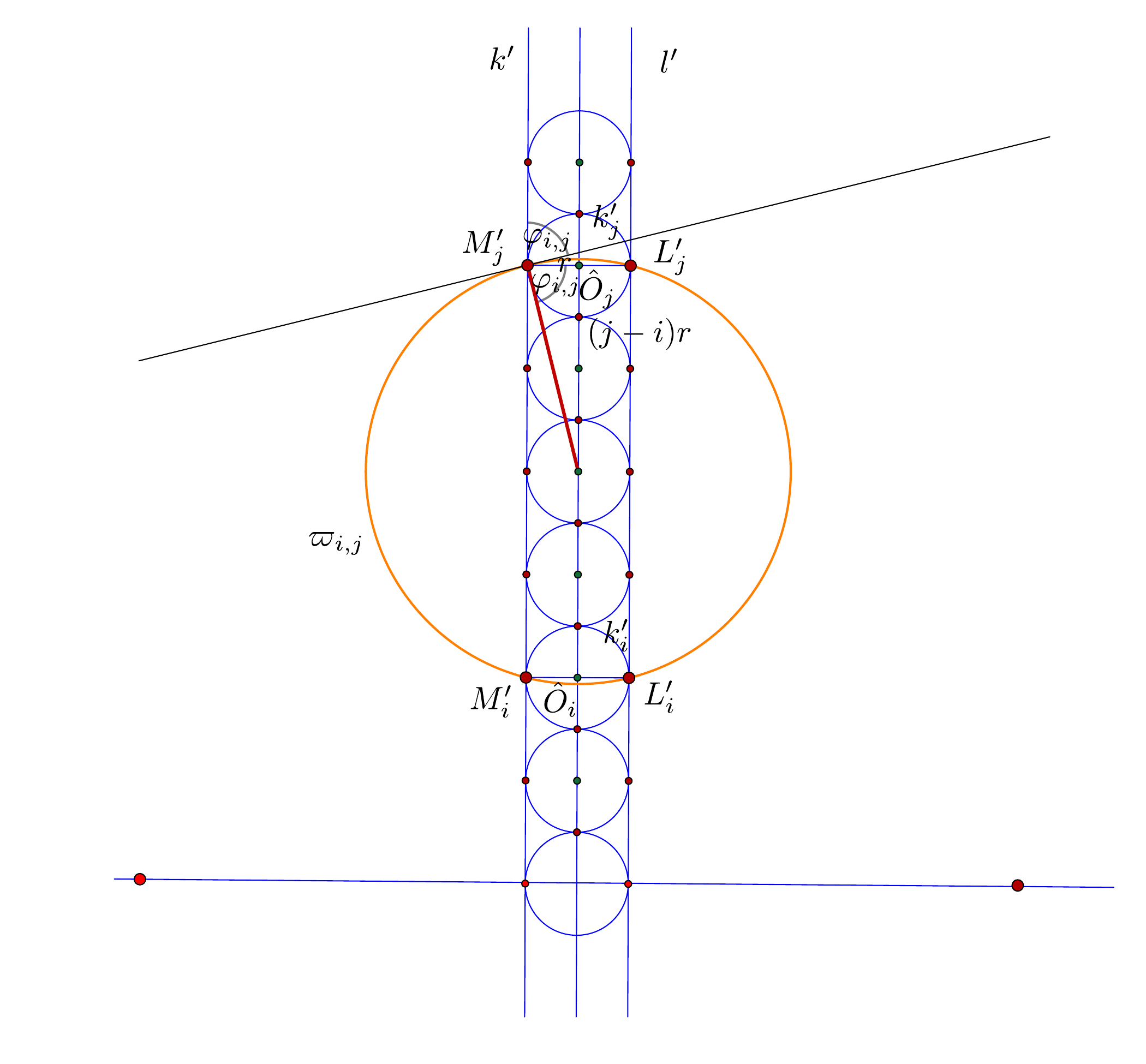}
    \caption{Proof of \ref{prop1}}\label{figl2}
\end{figure}

 The intersection angle of this circle $\varpi_{i,j}$ and $k'_j$ has $\tan (j-i)$, see Figure \ref{figl2} and the claim follows. \hfill $\square$

\begin{proposition} \label{prop2} For all $i, j\geq 0$ the circumscribed circle $\omega_{i, j}$ of  $\triangle L_i M_i N_j$ tangents both circles $k_j$ and $k_{j+1}$ at $N_j$.
\end{proposition}

\textsl{Proof:} From the previous proofs, $\triangle L'_i M'_i N'_j$ is isosceles, and its circumcenter lie on the line passing through the centers of all circles $k'_i$ in the chain, so it is tangent to $k'_j$ and $k'_{j+1}$ at $N'_j$. It yields the desired result. \hfill $\square$

The following theorem reveals exciting properties of the circles $\omega_{i, j}$ and $\varpi_{i, j}$ when $i-j=k$ where $k$ is a fixed integer.

\begin{figure}[!ht]
  \centering
    \includegraphics[width= 0.8 \textwidth]{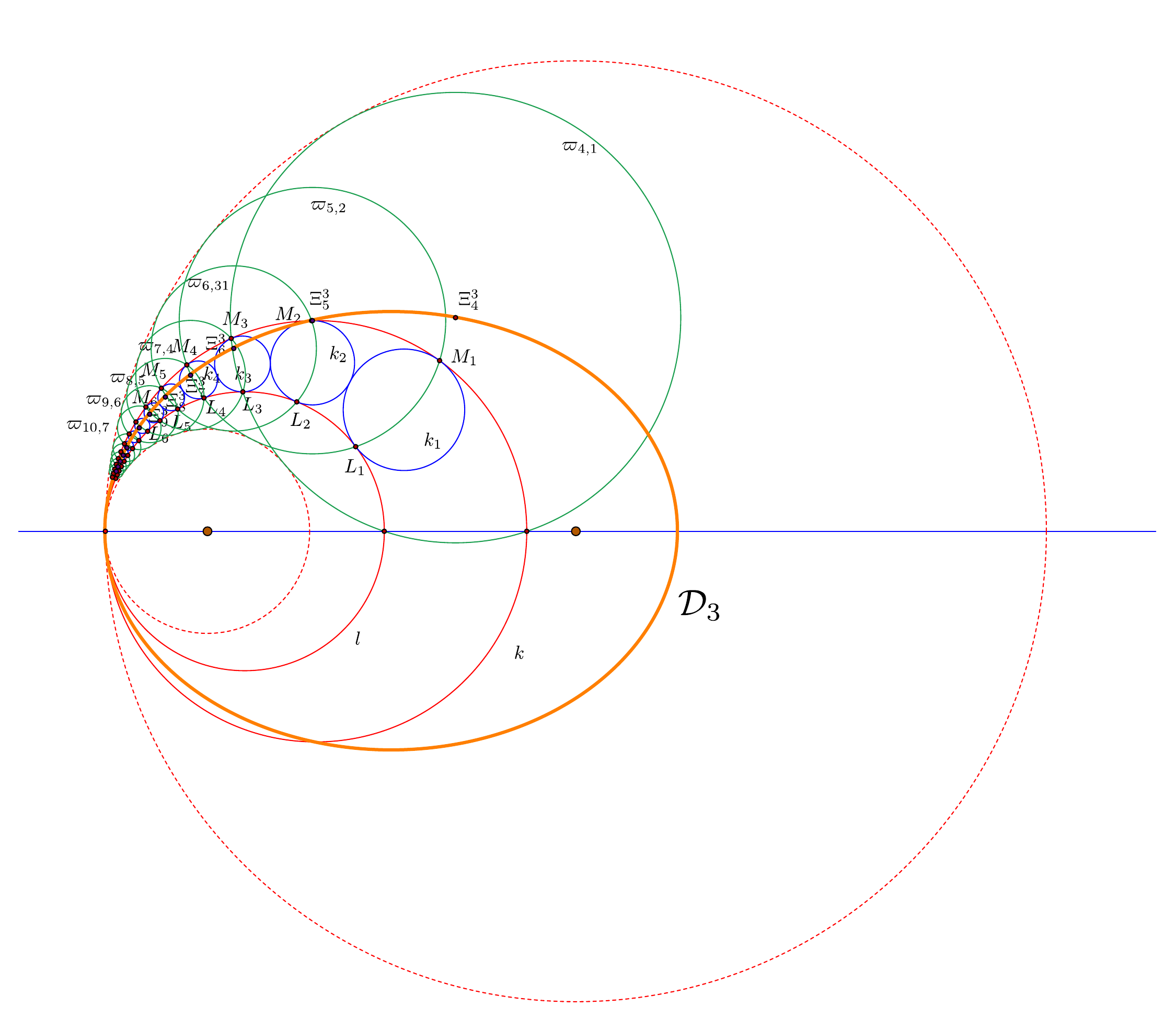}
    \caption{Theorem \ref{ceneters1}}\label{figl3}
\end{figure}

\begin{theorem}\label{ceneters1} \begin{itemize}
                  \item The centers $\Omega^k_i$ of all circles $\omega_{i, i-k}$ lie on the same ellipse $\mathcal{C}_k$ whose focuses lie on the line through the centers of $l$ and $m$.
                  \item The centers $\Xi^k_i$ of all circles $\varpi_{i, i-k}$ lie on the same ellipse $\mathcal{D}_k$ whose focuses lie on the line through the centers of $l$ and $m$, see Figure \ref{figl3}.
                  \item All ellipses $\mathcal{C}_k$ and $\mathcal{D}_k$  are mutually tangent and tangent to the circles $l$ and $m$ at $A$.
                \end{itemize}
\end{theorem}

\textsl{Proof:} We continue to use the same inversion $\psi$ from the previous proof.

 \begin{figure}[!ht]
  \centering
    \includegraphics[width= 0.45 \textwidth]{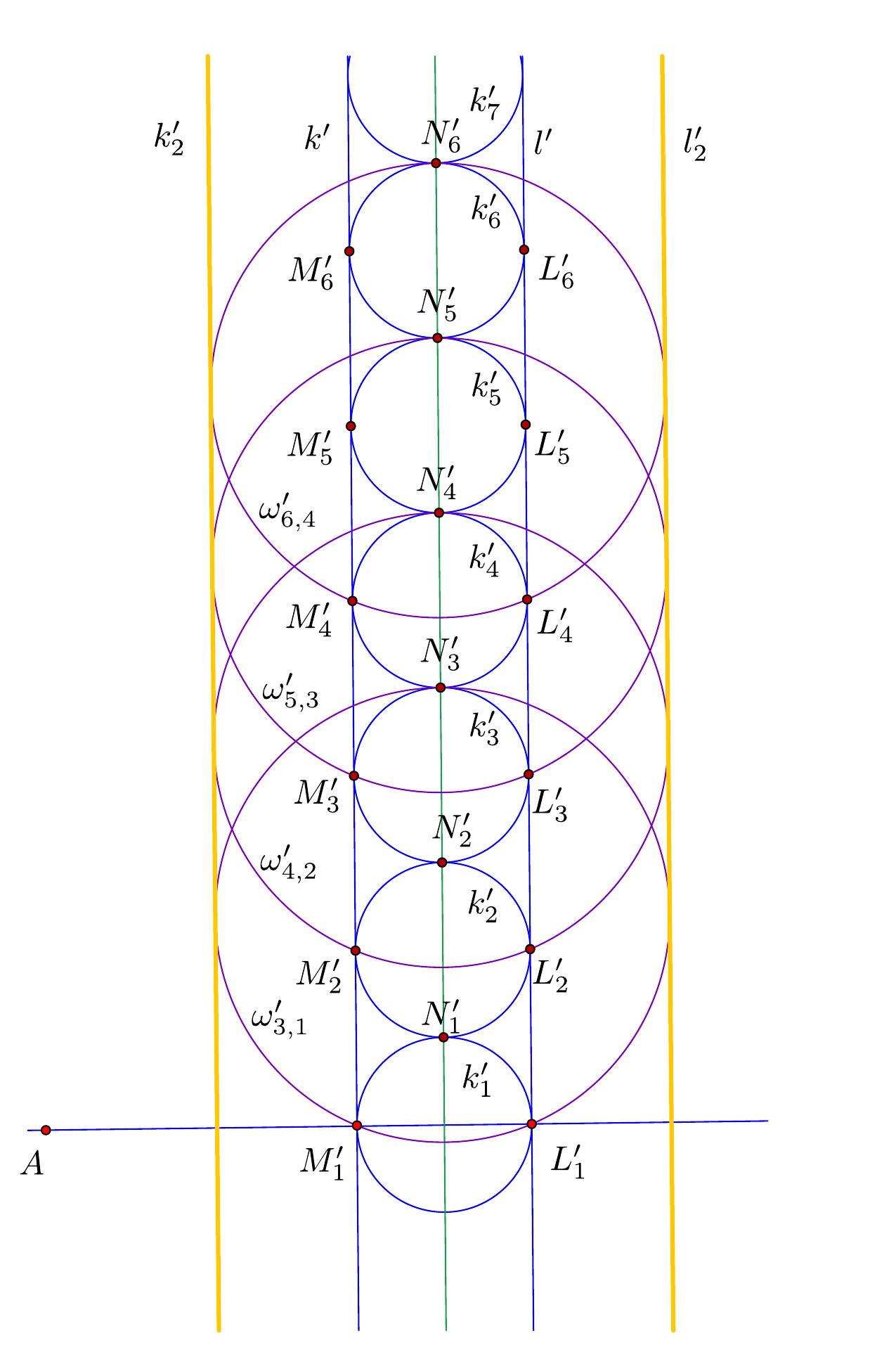}
    \caption{Proof of Theorem \ref{ceneters1}}\label{figl4}
\end{figure}

Observe (Figure \ref{figl4}) that there exists lines $l'_k$ and $m'_k$ parallel to $l'$ and $m'$, and tangent to all circles $\omega'_{i, i-k}$. It implies that all circles $\omega_{i, i-k}$ are tangent to the circles $l_k$ and $m_k$, the pre-images of $l'_k$ and $m'_k$, respectively. $l_k$ and $m_k$ are tangent to $l$ and $m$ at $A$, and their centers $U_k$ and $V_k$ lie on the line through the centers of $l$ and $m$. For each $i$, $\omega_{i, i-k}$ is externally tangent to $l$ and internally to $m$, so the sum of distances $\Omega^k_i U_k$ and $\Omega^k_i V_k$ is equal to the sum of the radiuses of $l_k$ and $m_k$. Therefore,  $\Omega^k_i$ lies on the same ellipse $C_k$ whose focuses are $U_k$ and $V_k$, such that it is tangent to $l$ and $m$ in $A$.

Analogously, we show the claim for the points $\Xi^k_i$.

The last part of the statement is a consequence of the fact that the point $A$, the focuses of all ellipses $\mathcal{C}_k$ and $\mathcal{D}_k$ are collinear. \hfill $\square$

 In \cite{Bar}, it is proved that for all $j >i \geq 0$, the points $L_j, L_i, M_i, M_j$ lie on the same circle $\varpi_{i, j}$ in the case of Steiner chains, and that they meet the circles $l$ and $m$ under the same angle. The following theorem is an analogous result of Theorem \ref{ceneters1} for the Steiner chain of circles.  Let $i-j=k$ where $k$ is a fixed integer.

\begin{figure}[!ht]
  \centering
    \includegraphics[width=0.8 \textwidth]{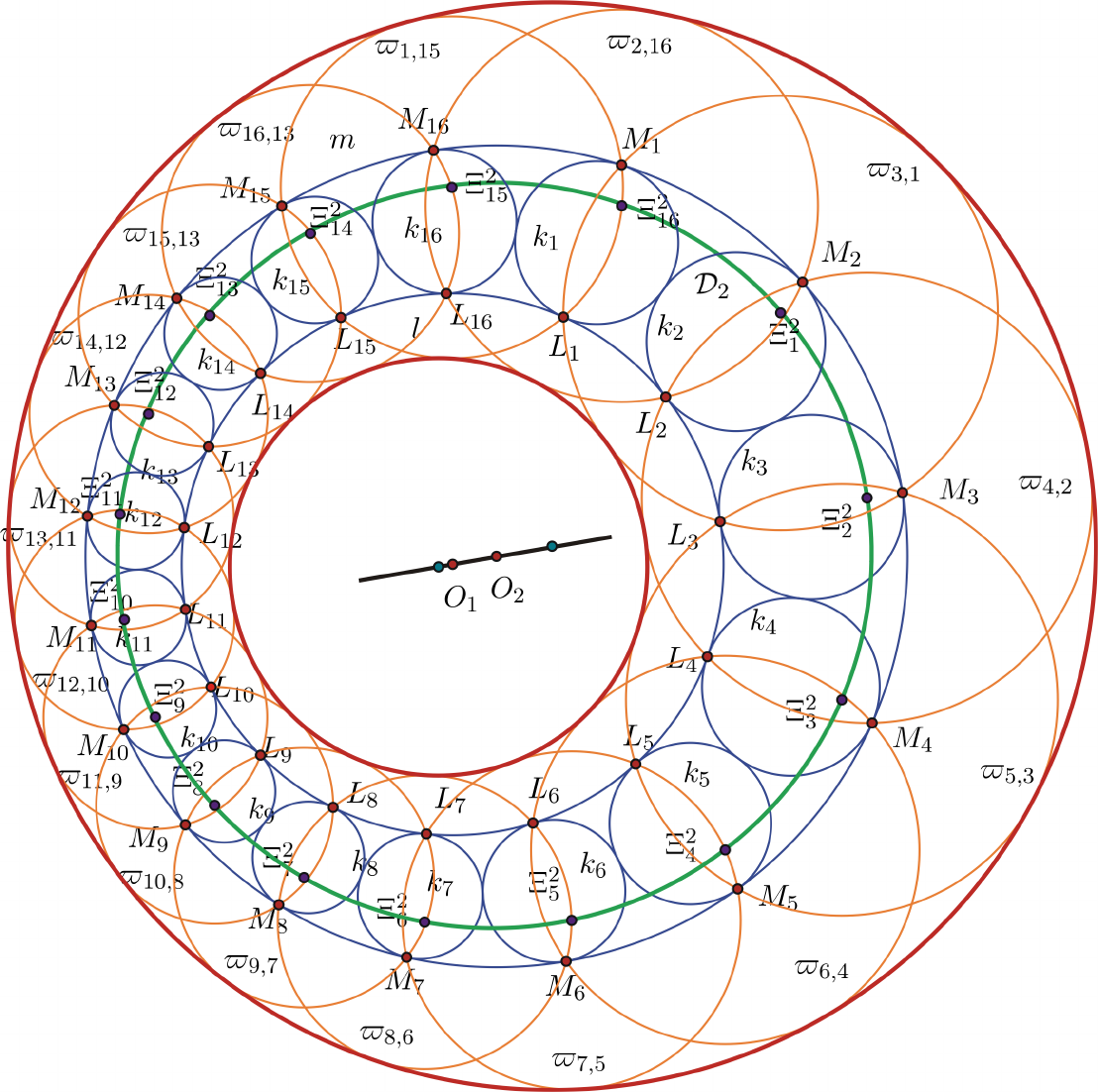}
    \caption{Theorem \ref{centers2} for $n=16$}\label{theorem27}
\end{figure}

 \begin{theorem} \label{centers2}In a Steiner chain, the centers $\Xi^k_i$ of all circles $\varpi_{i, i-k}$ lie on the same ellipse $\mathcal{D}_k$ whose focuses lie on the line through the centers of $l$ and $m$ (see Figure \ref{theorem27}).
 \end{theorem}

 \textsl{Proof:} We apply an inversion that maps $l$ and $m$ into two concentric circles. The existence of $\varpi_{i, i-k}$ has been proved in \cite[Proposition~3.1]{Bar}. Because of symmetry, the images of $\varpi_{i, i-k}$ are congruent and on the same distance from the common centre of $l'$ and $m'$. Therefore, there exist two concentric circles with $l^\prime$ and $m^\prime$ and tangent to all the images of $\varpi_{i, i-k}$. Consequently, two circles are tangent to $\varpi_{i, i-k}$, and their centres lie on the line through the centres of $l$ and $m$. Obviously, the centres of the two circles are focuses of the ellipse $\mathcal{D}_k$ passing through the centres of $\varpi_{i, i-k}$. \hfill $\square$

\section{Theorems on two orthogonal Pappus chains}

Let $k_1^1, \dots, k_n^1, \dots $ and $k_1^2, \dots, k_n^2, \dots $ be two Pappus chains of circles, touching $l_1$ and $m_1$, and $l_2$ and $m_2$, respectively.

\begin{definition}\label{def2} Two Pappus chains of circles of circles are called orthogonal if the circles $l_1$ and $m_1$ are orthogonal to the circles  $l_2$ and $m_2$, and $k_i^1=k_2^j$ for some $i, j$.
\end{definition}

The configuration of two orthogonal Pappus chains is special. Constraint on the circles $l_1$, $l_2$, $m_1$ and $m_2$ implies the following

\begin{proposition}
\begin{itemize}
                   \item The circles $l_1$, $l_2$, $m_1$ and $m_2$ pass through a common point $W$ that is the tangency point of $l_1$ and $m_1$, and $l_2$ and $m_2$,
                   \item The line through the centers of $l_1$ and $m_1$ is perpendicular to the line through the centers  $l_2$ and $m_2$ in $W$.
                 \end{itemize}
\end{proposition}

 \textsl{Proof:} The second claim is a corollary of the first, so it is enough to show only the first one. Let $W$ be the tangency point of $l_1$ and $m_1$. Let us apply an inversion with the centre in $W$. The circles $l_1$ and $m_1$ are mapped to the parallel lines $l'_1$ and $m'_1$. The image of $l_2$ is either a circle or a line orthogonal to $l'_1$ and $m'_1$. If the image is a circle $l'_2$, its centre lies on the line $l'_1$ since $l_2\perp l_1$. Similarly, the centre of $l'_2$ must lie on the line $m'_1$. But this is impossible since $l'_1\parallel m'_1$. Therefore, the image of $l_2$ is a line, so $l_2$ passes through $W$. Analogously, $m_2$ pass through $W$. \hfill $\square$

 The following property is a corollary of the above proposition.

 \begin{proposition} Let the circle $l_2$ intersects $l_1$ and $m_1$  in $W_1$ and $W_2$ distinct from $W$, and the circle $m_2$ intersects $l_1$ and $m_1$  in $W_3$ and $W_4$, respectively. Then the points $W_1$, $W_2$, $W_3$ and $W_4$ lie on the same circle meeting all the circles $l_1$, $l_2$, $m_1$ and $m_2$ under angle $45^\circ$.
 \end{proposition}

  \textsl{Proof:} After applying an inversion with the centre in $W$, the points $W_1$, $W_2$, $W_3$ and $W_4$ are mapped into the vertices of a rectangle, so they belong to a circle. Since two chains are orthogonal, there is a circle inscribed in the rectangle $W'_1 W'_2 W'_3 W'_4$ so $W'_1 W'_2 W'_3 W'_4$ is square, and the claim follows. \hfill $\square$

From the next lemma, it follows that any circle in a Pappus chain is the common circle for two orthogonal Pappus chains.

\begin{lemma} Let $k_1^1, \dots, k_n^1, \dots $ be a Pappus chain. For every positive integer $i$ there is an orthogonal Pappus chain $k_1^2, \dots, k_n^2, \dots $ such that $k_1^2=k_i^1$
\end{lemma}

  \textsl{Proof:} We apply an inversion with the centre in $W$. The chain $k_1^1, \dots, k_n^1, \dots $ is mapped to the chain of equal circles filling the strip between the lines $l'_1$ and $m'_1$. Let us construct two tangents  $l'_2$ and $m'_2$, of the ${k'}_i^1$ perpendicular to $l'_1$ and $m'_1$. The strip between $l'_2$ and $m'_2$ may be also filled with the congruent circles ${k'}_1^2={k'}_i^1, \dots, {k'}_n^2, \dots $. However, the image of this chain of circles and $l'_2$ and $m'_2$ is the wanted orthogonal chain. \hfill $\square$

Indeed, for orthogonal Pappus chains without loss of generality, we may assume that they have the first circles in the chains. In the rest of the paper, we assume that unless it is stated otherwise.

\begin{figure}[!ht]
  \centering
    \includegraphics[width=0.8 \textwidth]{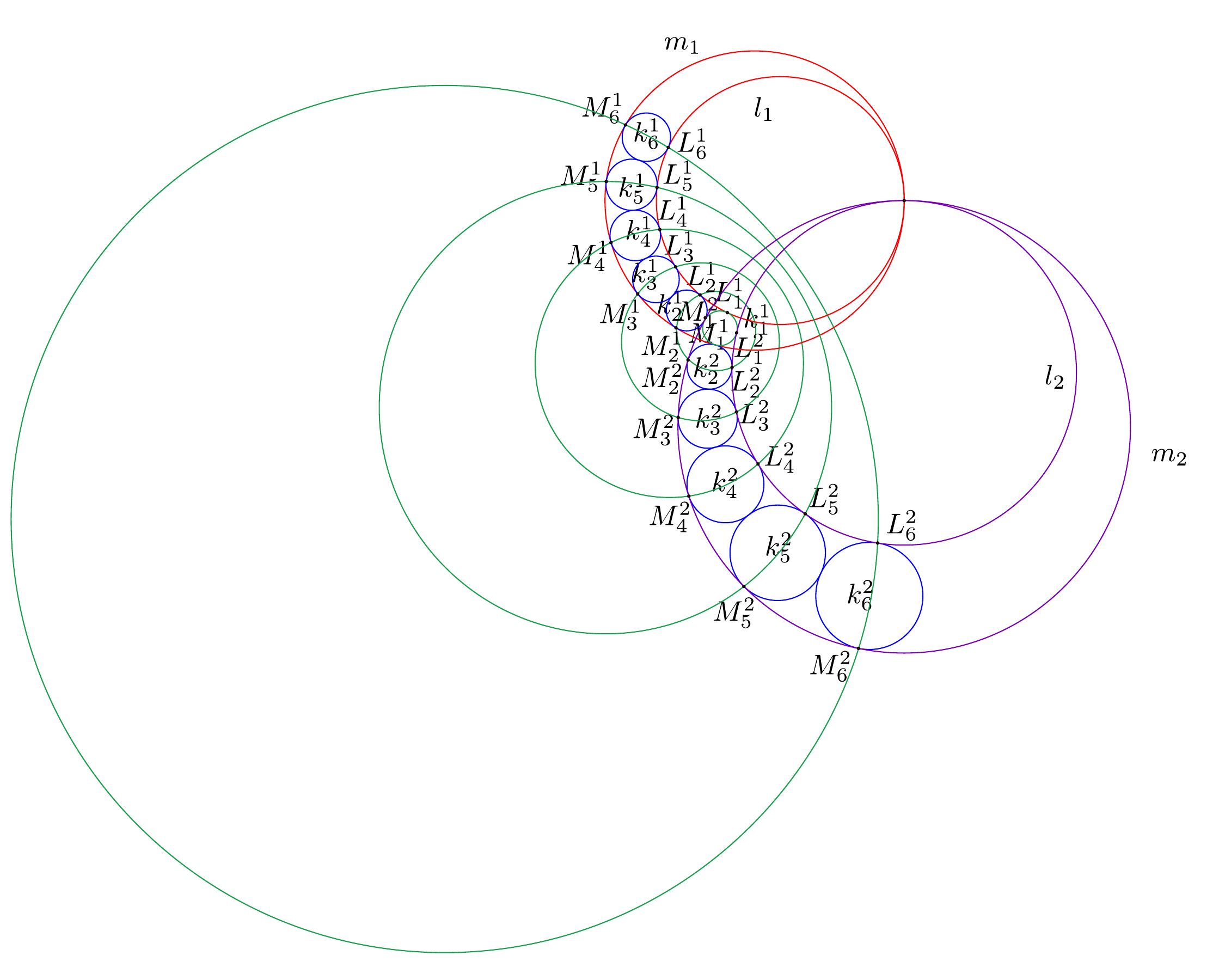}
    \caption{Theorem \ref{big1}}\label{t36a}
\end{figure}

\begin{theorem}\label{big1} Let $k_1^1, \dots, k_n^1, \dots $ and $k_1^2, \dots, k_n^2, \dots $ be two orthogonal Pappus chains of circles. Then:

\begin{itemize}
  \item For each $i$, the points $L_i^1, M_i^1, L_i^2, M_i^2$ lie on the same circle $c_i$. There exists a circle (or a line) passing through $W$ orthogonal to all circles $c_i$ and $k_1^1$. The centres of circles ll circles $c_i$ and $k_1^1$ lie on the same line; see Figure \ref{t36a}.
  \item For all $i, j$ the points $L_i^1, L_j^1, L_i^2, L_j^2$ lie on the same circle $c_{i, j}$,
  \item For all $i, j$ the points $M_i^1, M_j^1, M_i^2, M_j^2$ lie on the same circle $\bar{c}_{i, j}$,
  \item The centers of all circles $c_{1, k}$, $c_{2, 2k+1}$, $c_{3, 3k+2}, \dots , c_{n, nk+n-1}, \dots $ lie on the same hyperbola $\mathcal{H}_k$ and the centers of  $\bar{c}_{1, k}$, $\bar{c}_{2, 2k+1}$, $\bar{c}_{3, 3k+2}, \dots , \bar{c}_{n, nk+n-1}, \dots $ lie on the same hyperbola $\mathcal{\bar{H}}_k$
\end{itemize}
\end{theorem}

 \textsl{Proof:} We apply an inversion with center in $W$. All circles from the chains are mapped into congruent circles distributed between two orthogonal strips. As a corollary, we have that the vertices of quadrilateral $L_i^1 M_i^1 M_i^2 L_i^2 $ are mapped to the vertices of isosceles trapezium ${L'}_i^1 {M'}_i^1 {M'}_i^2 {L'}_i^2 $. Isosceles trapezium has a circumscribed circle, so the points $L_i^1, M_i^1, L_i^2, M_i^2$ lie on the same circle $c_i$. All segments ${L'}_i^1 {L'}_i^2$ and ${M'}_i^1 {M'}_i^2$ have a common perpendicular bisector that passes through the centre of ${k'}_1^1$ and this line is mapped to circle or line passing through $W$ and orthogonal to $k_1^1$ and all circles $c_i$. Observe that the circles $c'_i$ are concentric, so the centres of $c_i$ lie on the same line.

 In the same fashion, ${L'}_i^1 {L'}_j^1 {L'}_j^2 {L'}_i^2$ and ${M'}_i^1 {M'}_j^1 {M'}_j^2 {M'}_i^2$ are isosceles trapeziums so the points $L_i^1, L_j^1, L_i^2, L_j^2$ are concyclic as well as $M_i^1, M_j^1, M_i^2, M_j^2$.

   \begin{figure}[!ht]
  \centering
    \includegraphics[width=0.7 \textwidth]{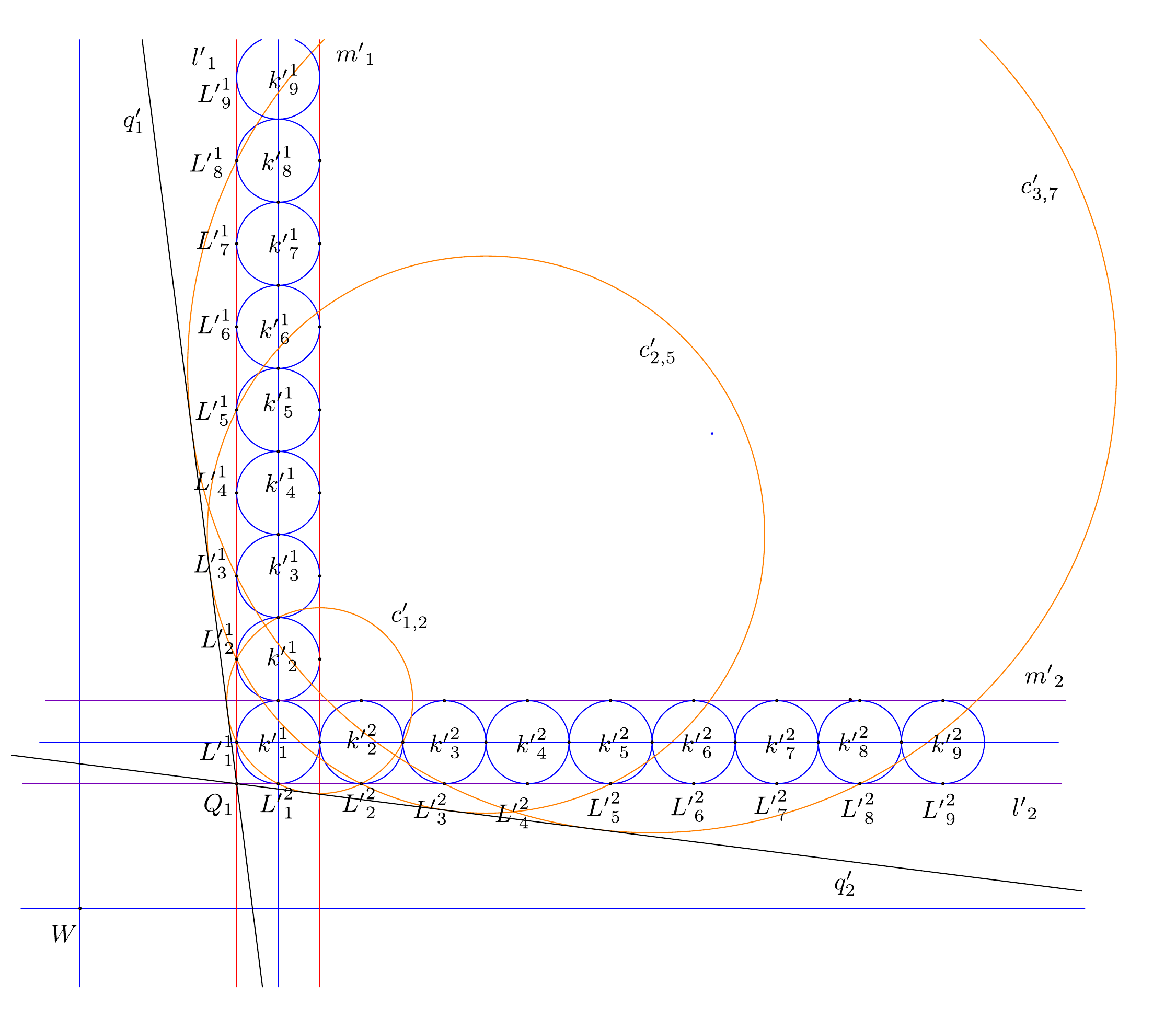}
    \caption{Proof of Theorem \ref{big1} for $k=2$}\label{t36b}
\end{figure}

 Let us observe that there exists a homothety with centre in the intersection point $Q_l$ of the lines $l'_1$ and $l'_2$ that sends ${L'}_1^1 {L'}_k^1 {L'}_k^2 {L'}_1^2$ to ${L'}_j^1 {L'}_{(2j-1)k}^1 {L'}_{(2j-1)k}^2 {L'}_j^2$ for all $j$, see Figure \ref{t36b}. Therefore, the tangents from $Q_1$ to $c'_{1, k}$ are also tangent to $c'_{j, j(k+1)-1}$ for all $j$. Indeed, all circles $c'_{1, k}$, $c'_{2, 2k+1}$, $c'_{3, 3k+2}, \dots , c'_{n, nk+n-1}, \dots $ have two common tangents $q'_1$ and $q'_2$ passing through $Q_1$, so the circles $c_{1, k}$, $c_{2, 2k+1}$, $c_{3, 3k+2}, \dots , c_{n, nk+n-1}, \dots $ are internally tangent to two circles $q_1$ and $q_2$ passing through the intersection points of circles $l_1$ and $l_2$. But, it implies that their centres lie on the same hyperbola $\mathcal{H}_k$ whose focuses are the centres of the circles $q_1$ and $q_2$.

We prove the claim for the centers of  $\bar{c}_{1, k}$, $\bar{c}_{2, 2k+1}$, $\bar{c}_{3, 3k+2}, \dots , \bar{c}_{n, nk+n-1}, \dots $ analogously. \hfill $\square$

\begin{figure}[!ht]
  \centering
    \includegraphics[width=0.8 \textwidth]{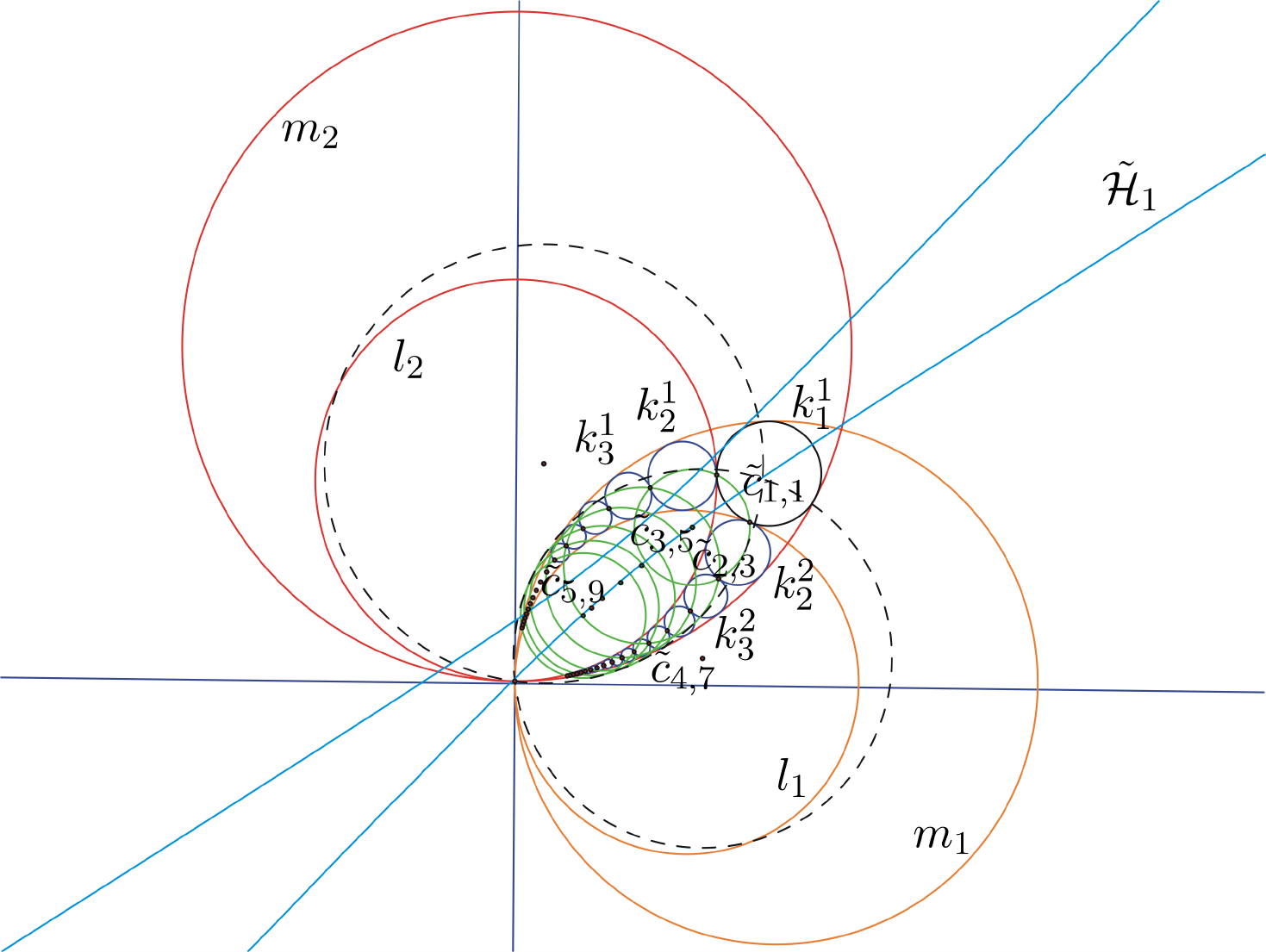}
    \caption{Theorem \ref{big2} for $k=2$}\label{mainteo}
\end{figure}

\begin{theorem}\label{big2} Let $k_1^1, \dots, k_n^1, \dots $ and $k_1^2, \dots, k_n^2, \dots $ be two orthogonal Pappus chains of circles. Then:

\begin{itemize}

  \item For all $i, j$ the points $N_i^1, N_j^1, N_i^2, N_j^2$ lie on the same circle $\tilde{c}_{i, j}$,
  \item The centers of all circles $\tilde{c}_{1, k}$, $\tilde{c}_{2, 2k+1}$, $\tilde{c}_{3, 3k+2}, \dots , \tilde{c}_{n, (nk+n-1}, \dots $ lie on the same hyperbola $\mathcal{\tilde{H}}_k$, see Figure \ref{mainteo}.
\end{itemize}
\end{theorem}

 \textsl{Proof:} Like the previous proof, we apply the inversion with the centre in $W$. Under the inversion,  the points $N_i^1, N_j^1, N_i^2, N_j^2$ are mapped to the vertices of isosceles trapezium ${N'}_i^1 {N'}_j^1 {N'}_j^2 {N'}_i^2$ for all $i, j$ so they are concyclic.

  \begin{figure}[!ht]
  \centering
    \includegraphics[width=0.7 \textwidth]{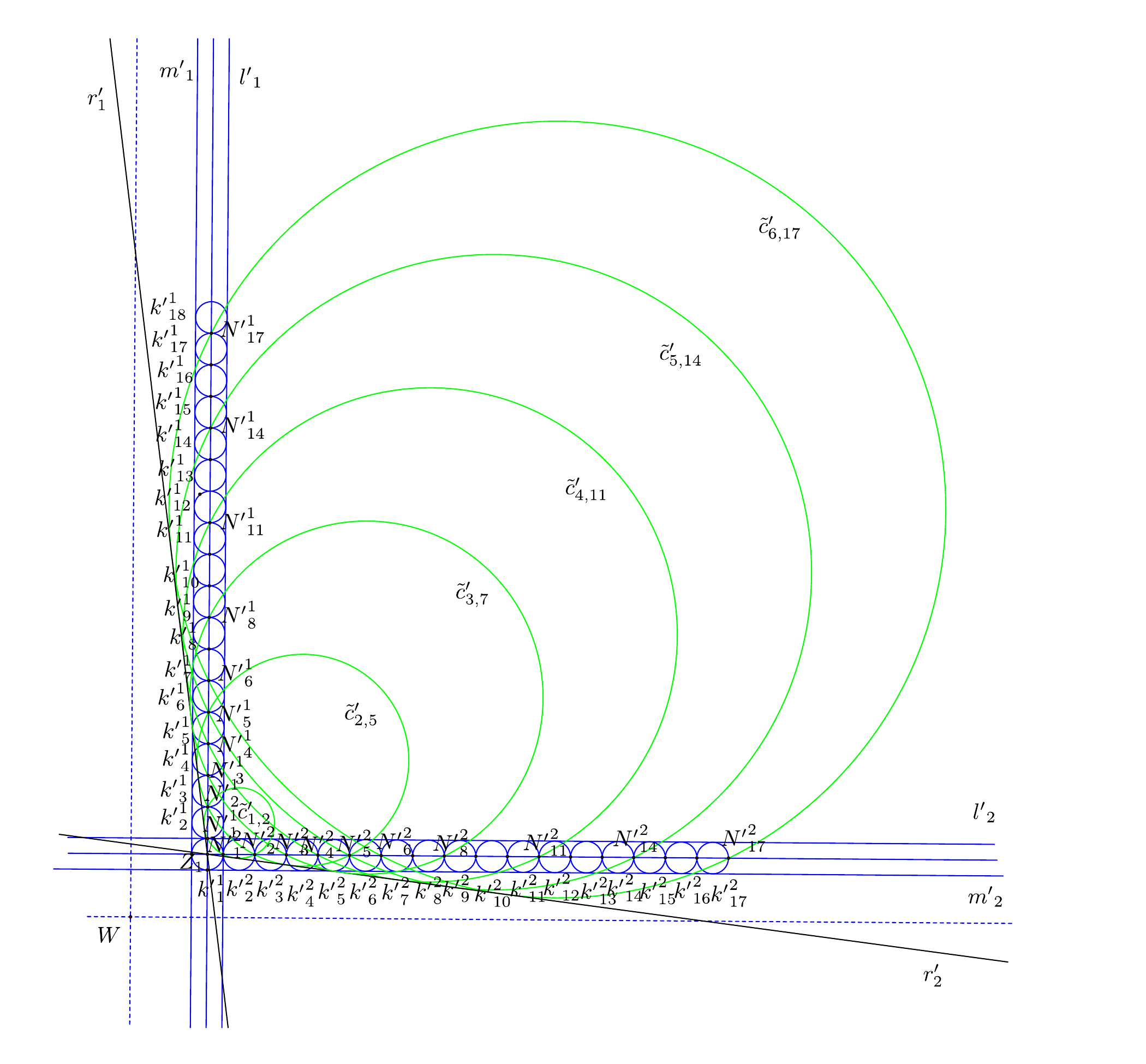}
    \caption{Proof of Theorem \ref{big2} for $k=2$}\label{mainteop}
\end{figure}

 Let $Z_1$ be the centre of the circle ${k'}_1^1$. For all $j, k$, there is a homothety that maps isosceles trapezium  ${N'}_1^1 {N'}_k^1 {N'}_k^2 {N'}_1^2$ to isosceles trapezium ${N'}_j^1 {N'}_{(2j-1)k}^1 {N'}_{(2j-1)k}^2 {N'}_j^2$. The latter observation implies that $\tilde{c}'_{1, k}$ is mapped to $\tilde{c}'_{j, j(k+1)-1}$ under the homothety, see Figure \ref{mainteop}. Therefore, all circles $\tilde{c}'_{1, k}$, $\tilde{c}'_{2, 2k+1}$, $\tilde{c}'_{3, 3k+2}, \dots , \tilde{c}'_{n, nk+n-1}, \dots $ have two common tangents $r'_1$ and $r'_2$ passing through $Z_1$, so the circles $\tilde{c}_{1, k}$, $\tilde{c}_{2, 2k+1}$, $\tilde{c}_{3, 3k+2}, \dots , \tilde{c}_{n, nk+n-1}, \dots $ are internally tangent to two circles $r_1$ and $r_2$ passing through the intersection points of circles $l_1$ and $l_2$. But, it implies that their centres lie on the same hyperbola $\mathcal{\mathcal{H}}_k$ whose focuses are the centres of the circles $r_1$ and $r_2$. The circles $r_1$ and $r_2$ are orthogonal to the circle $k_1^1$. \hfill $\square$

\section{Nonexistence of two orthogonal Steiner chains  and orthogonal Steiner and Pappus chains}

The notion of orthogonality introduced in Definition \ref{def2} is too strong to be generalized to orthogonality: two  Steiner chains or a Steiner chain and a Pappus chain. It turns out that the request to have orthogonal outer and inner circles in two chains is impossible unless in the case of two Pappus chains having the common tangent point of the inner and the outer circles. This observation follows directly from the next lemma.

\begin{lemma}Let $l$ and $k$ be circles such that circle $l$ lies inside $k$. Let $m_1$ and $m_2$ be two circles orthogonal to both $l$ and $k$. Then $m_1$ and $m_2$ have two intersection points.
\end{lemma}

  \textsl{Proof:} We consider inversion that maps $l$ and $k$ into two concentric circles $l'$ and $k'$ with center in $X'$. Let $m'_1$ be the image of $m_1$. If $m'_1$ is a circle with center in $Y'_1$ then ${X' Y'_1}^2=r_{l'}^2+r_{m'_1}^2=r_{k'}^2+r_{m'_2}^2$ where $r_{l'}$, $r_{k'}$ and $r_{m'_1}$ are the radiuses of $l'$, $k'$ and $m'_1$. But, the equalities cannot hold simultaneously since  $r_{l'}\neq r_{k'}$. Thus, $m'_1$ is a line passing through $X'$ because of orthogonality. Analogously, $m'_2$ is a line passing through $X'$.

  Therefore, $m_1$ and $m_2$ pass through the centre of the inversion and the preimage of $X'$. The claim is proved.
  \hfill $\square$

\begin{figure}[!ht]
  \centering
    \includegraphics[width=0.55 \textwidth]{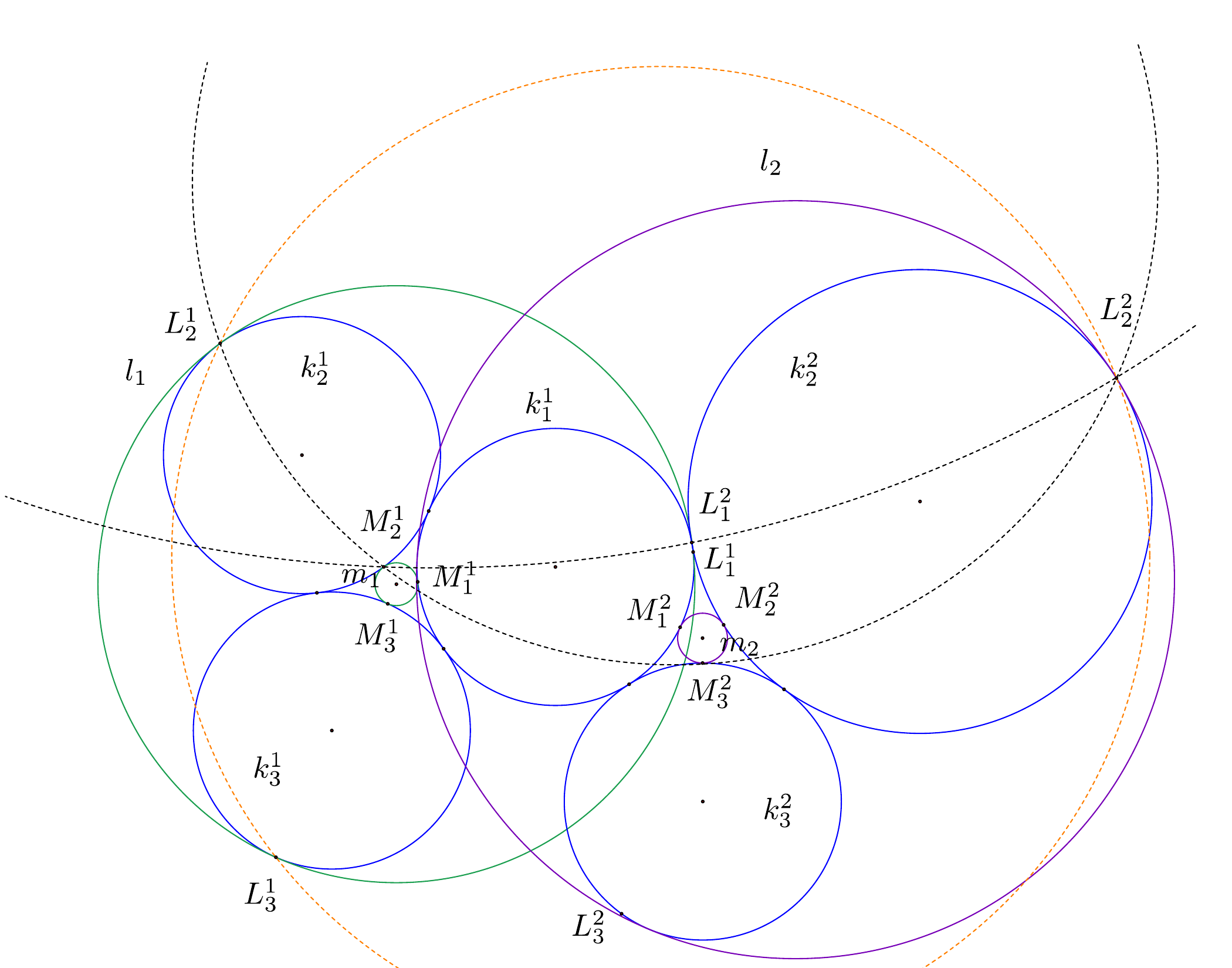}
    \caption{Counterexample for Theorems \ref{big1} and \ref{big2} for two Steiner chains with common circle}\label{kon}
\end{figure}

A natural question arising here is whether we can find analogues of Theorems \ref{big1} and \ref{big2} in some cases of Steiner chains having a common circle. Do they exist, and if so, under which constraints on Steiner chains do such constructions exist? An easy experiment in Cinderella software shown in Figure \ref{kon} says the answer is negative in the general case.

The question is, of course, far more complex than in the case of orthogonal Pappus chains, as there is no inversion sending the four fixed circles $l_1$, $m_1$, $l_2$ and $m_2$ to the lines. However, the classical proof of Steiner porism suggests that a hypothetical desired configuration would be an inversion that maps circles so that  $l'_1$ and $m'_1$, as well as $l'_2$ and $m'_2$ are concentric, while the circles in both chains are congruent. This situation is, of course, possible, but in particular cases.

With the above idea in mind, it is not hard to find configurations of two Steiner chains having a common circle in the chains satisfying analogues of Theorems \ref{big1} and \ref{big2}, as confirmed by the experiments done in Cinderella.

It is not hard to prove the following claim.

\begin{theorem}\label{zad}
Let $k_1^1, \dots, k_n^1 $ and $k_1^2, \dots, k_n^2 $ be two Steiner chains with common circle $k_1^1=k_2^1$, such that there exists an inversion  $\psi$ mapping $l_1$ and $m_1$, and $l_2$ and $m_2$ in two concentric circles. Then:

\begin{itemize}

  \item For all $i, j$ the points $N_i^1, N_j^1, N_i^2, N_j^2$ lie on the same circle $c_{i, j}$,

  \item For all $i, j$ the points $M_i^1, M_j^1, M_i^2, M_j^2$ lie on the same circle $c^M_{i, j}$,

  \item For all $i, j$ the points $L_i^1, L_j^1, L_i^2, L_j^2$ lie on the same circle $c^L_{i, j}$,

  \item For all $i$ the points$L_i^1, M_i^1, L_i^2, M_i^2$ lie on the same circle $c_i$.

\end{itemize}

\end{theorem}

\begin{figure}[!ht]
  \centering
    \includegraphics[width=0.75 \textwidth]{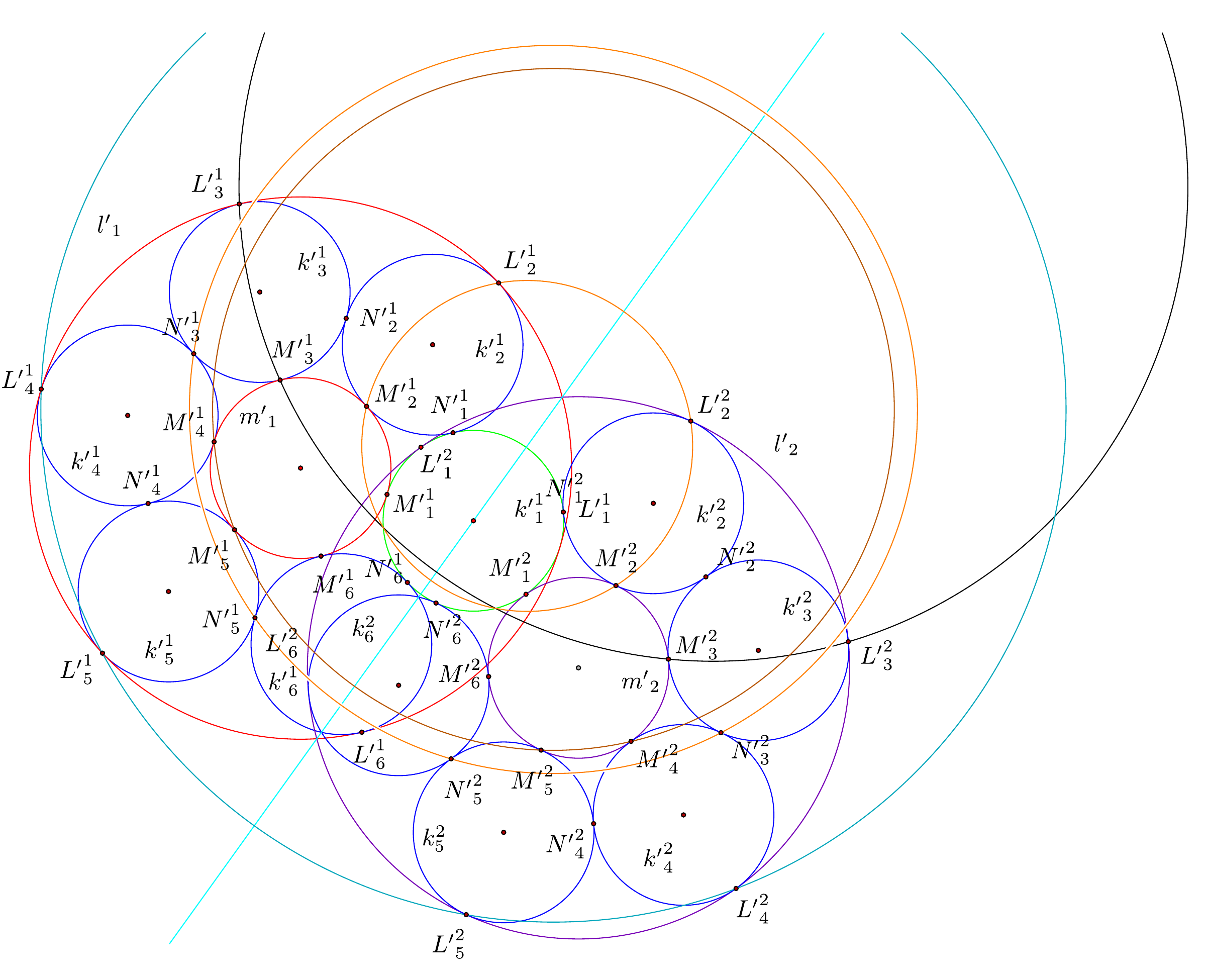}
    \caption{Proof of Theorem \ref{zad} }\label{pri}
\end{figure}

\textsl{Proof:} Consider the image under inversion $\psi$, Figure \ref{pri}.

The first Steiner chain can mapped to the second so that  $l'_1$ and $m'_1$ are mapped $l'_2$ and $m'_2$, and ${k'}_i^1$ in ${k'}_{n+2-i}^2$ for all $2\leq i \leq n$, by some rotation around the centre of ${k'}_1^1$. Composing the rotation with the symmetry with respect to the line through the centers ${k'}_1^1$ and ${l'}_2$, we deduce that there is a symmetry with respect to a line through the centre of circle  ${k'}_1^1$ that maps ${k'}_i^1$ into ${k'}_i^2$ for every $i$.

Therefore, the quadrilaterals ${N'}_i^1 {N'}_j^1 {N'}_i^2 {N'}_j^2$, ${L'}_i^1 {L'}_j^1 {L'}_i^2 {L'}_j^2$, ${M'}_i^1 {M'}_j^1 {M'}_i^2 {M'}_j^2$ and ${L'}_i^1 {M'}_i^1 {L'}_i^2 {M'}_i^2$ are isosceles trapeziums, and the claim follows. \hfill $\square$

Experimenting with Cinderella, we discovered the following exciting result.

\begin{figure}[!ht]
  \centering
    \includegraphics[width=0.55 \textwidth]{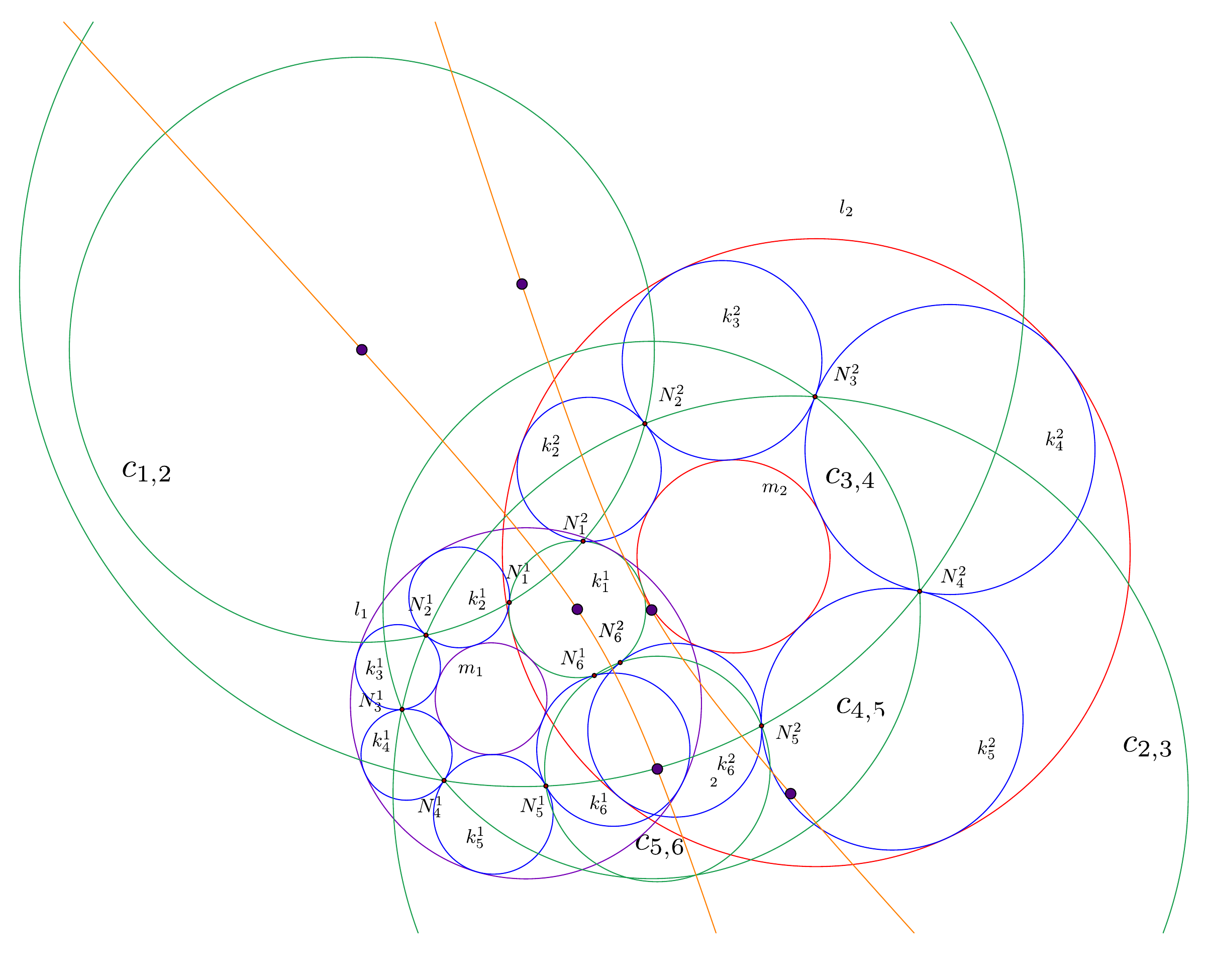}
    \caption{Theorem \ref{hip}}\label{hipo}
\end{figure}

\begin{theorem}\label{hip}
Let $k_1^1, \dots, k_n^1 $ and $k_1^2, \dots, k_n^2 $ be two Steiner chains with common circle  $k_1^1=k_2^1$ such that there is an inversion $\psi$ sending $l_1$ and $m_1$, as well as $l_2$ and $m_2$ in two concentric circles. Then:

\begin{itemize}

  \item Centers of the circles  $c_{1, k}$, $c_{2, k+1}$, $c_{3, k+2}, \dots , c_{n, k-1}$ lie on the same conic $\mathcal{H}_k$, see Figure \ref{hipo}

  \item Centers of the circles  $c^L_{1, k}$, $c^L_{2, k+1}$, $c^L_{3, k+2}, \dots , c^L_{n, k-1}$ lie on the same conic $\mathcal{H}^L_k$,

  \item Centers of the circles  $c^M_{1, k}$, $c^M_{2, k+1}$, $c^M_{3, k+2}, \dots , c^M_{n, k-1}$ lie on the same conic $\mathcal{H}^M_k$,

  \item Centers of the circles  $c_1$, $c_2$, $c_3, \dots , c_n$ lie on the same conic $\mathcal{C}$.
\end{itemize}

\end{theorem}

\textsl{Proof:} In the inversion $\psi$, the images of the circles  $c^L_{1, k}$, $c^L_{2, k+1}$, $c^L_{3, k+2}, \dots , c^L_{n, k-1}$ are the circles whose centres lie on the line of symmetry through the centre of the common circle of two chains described in the previous proof. The image circles have equal common  chords with $l'_1$ and $l'_2$. The claim now follows from Lemma \ref{najl}, which will be proved next.

Lemma \ref{najl} implies other results too. \hfill $\square$

Now we formulate and prove Lemma \ref{najl}.

\begin{lemma}\label{najl} Let $k$ and $\Omega$ be two given circles,  $p$ be a line not passing through the center of $k$ and $d$ given length. For a chord of the length $d$ of $k$, consider the circle with centre on the line $p$ passing through the endpoints of the chords and its image under the inversion with respect to $\Omega$. The locus of the centre of the image circle is a conic.
\end{lemma}

\textsl{Proof:} Without lose of generality, we can assume that $k$ is the unit circle, and the equation of line $p$ in the complex plane is $z+\bar{z}=2 a$, where $a\geq 0$. Let the affix of the centre of $\Omega$ be $c=p+i q$ and the radius of $\Omega$ be $\rho$. For a point $M$ on the line $p$, there are two circles centred in $M$ having the common chord of length $d$ with $k$. Let their radiuses be $r_1$, and Under the inversion with respect to $\Omega$, they are, in the general case, mapped to two circles whose centres have affixes $c+\rho^2 \frac{m-c}{|m-c|^2-r_1^2}$ and $c+\rho^2 \frac{m-c}{|m-c|^2-r_2^2}$, where $m=a+i y$ is the affix of the point $M$.

We have that $r_1^2, r_2^2= \frac{d^2}{4}+\left(|m|\mp\sqrt{1-\frac{d^2}{4}}\right)^2=|m|^2\pm 2 \sqrt{1-\frac{d^2}{4}} |m|+1$. Therefore, the loci of the centres are the points whose affixes satisfy  $c+\rho^2 \frac{m-c}{-m\bar{c}-\bar{m}c+|c|^2\pm 2 \sqrt{1-\frac{d^2}{4}} |m|-1}$. Adding the vector $c$ is translation, multiplication by $\rho$ is homothety, so it is enough to show that the set of points $\left(\frac{a-p}{Z}, \frac{y-q}{Z}\right)$ is a conic, where $Z=-2ap-2qy+p^2+q^2\pm 2 \sqrt{1-\frac{d^2}{4}} \sqrt{a^2+y^2}-1$.

Let $X=\frac{a-p}{Z}$ and $Y=\frac{y-q}{Z}$. Then $$y=q+(a-p)\frac{Y}{X}.$$ Also, $Z=\frac{a-p}{X}$. Put new constants $w=p^2+q^2-2ap-1$ and $s= 2 \sqrt{1-\frac{d^2}{4}}$ so $Z= w-2qy\pm s\sqrt{a^2+y^2}$. Now, we substitute $y$ and $Z$ to get $$\frac{2q(a-p)Y+(2q^2-w) X+a-p}{X}=\pm s \sqrt{\frac{2 a^2 X^2+2 a(a-p) X Y+(a-p)^2 Y^2}{X^2}}.$$ After squaring and multiplying by $X^2$ we see obtain the equation of a conic. \hfill $\square$

The results in this section trigger the following question:

\begin{problem} Find necessary conditions when two Steiner chains having a common circle in the chains admit higher order incidence results like in Theorem \ref{hip}?
\end{problem}

\begin{center}\textmd{\textbf{Acknowledgements} }
\end{center}

\medskip
The first author was partially supported by the Ministry of Science, Technological Development and innovations of the Republic of Serbia through the institutional funding of the Mathematical Institute SANU.

\begin{center}\textmd{\textbf{Data availability} }
\end{center}

\medskip
Data sharing not applicable to this article as no datasets were generated or analyzed during the current study.

\begin{center}\textmd{\textbf{Code availability} }
\end{center}

\medskip
Not applicable.

\begin{center}\textmd{\textbf{Conflict of interest} }
\end{center}

\medskip
The authors have no relevant financial or non-financial interests to disclose.

\end{document}